\documentclass[conference,a4paper%
]{IEEEtran}

\usepackage[utf8]{inputenc} 
\usepackage[T1]{fontenc}
\usepackage{url}
\usepackage{ifthen}
\usepackage{cite}
\usepackage[cmex10]{amsmath} 
\usepackage{amsmath}
\usepackage{amssymb}
\usepackage{amsthm}
\usepackage[para]{threeparttable}
\usepackage{general-private}
\usepackage{special-private}
\usepackage{thm-config}

\begin{document}
\title{Exact Calculation of Normalized Maximum Likelihood Code Length Using Fourier Analysis} 

\author{%
  \IEEEauthorblockN{Atsushi Suzuki and Kenji Yamanishi}
  \IEEEauthorblockA{The University of Tokyo\\
                    Graduate School of Information Science and Technology\\
                    Bunkyo, Tokyo, Japan\\
                    Email: atsushi.suzuki.rd@gmail.com, yamanishi@mist.i.u-tokyo.ac.jp}
}

\maketitle

\begin{abstract}
  The normalized maximum likelihood code length has been widely used in model selection, and its favorable properties, such as its consistency and the upper bound of its statistical risk, have been demonstrated. This paper proposes a novel methodology for calculating the normalized maximum likelihood code length on the basis of Fourier analysis. Our methodology provides an efficient non-asymptotic calculation formula for exponential family models and an asymptotic calculation formula for general parametric models with a weaker assumption compared to that in previous work.
\end{abstract}


\section{Introduction}
\subsection{Background and Our Contribution}
The normalized maximum likelihood code length (NML code length) is an extension of self-entropy, in which a set of distributions is given instead of the true distribution.
When the true distribution is known, the lower bound of the mean length of codings for a random variable is given by the Shannon entropy of its probability distribution, and the lower bound is attained by the self-entropy \cite{Shannon:1948}.
This optimal code or the self-entropy is also interpreted as the solution of a (trivial) optimization problem of a log redundancy with respect to the code $\CodeLength{}{}$, with the Kraft-McMillan inequality \cite{Kraft:1949} \cite{McMillan:1956} as a constraint, as follows:

\begin{equation}
\EqnLabel{ShannonMinMax}
\begin{split}
& \min_{\CodeLength{}{}} \max_{\DataSequence} \Bracket{\CodeLength{\DataSequence}{} - \Paren{- \log \TrueDensityFunc{\DataSequence}}} \\
& \ST \int \DataDiff \Exponential{- \CodeLength{\DataSequence}{}} \le 1,
\end{split}
\end{equation}
where $\DataSequence \DefEq \Range{\DataVec_{1}}{\DataVec_{2}}{\DataVec_{\DataCnt}}$ is a data sequence and $\TrueDensityFunc{}$ denotes the probability density function of the data-generating distribution. Apparently, the optimum code length is given by $\CodeLength{\DataSequence}{} = - \log \TrueDensityFunc{\DataSequence}$ (self-entropy). Note that we discuss cases of continuous random variables in this paper. Further, the base of the logarithm is $\mathrm{e}$, and the natural unit of information is used in this paper.

The optimization problem \EqnRef{ShannonMinMax} or the original Shannon entropy deals with the case in which the true distribution is known. When a set of distributions $\DensityDomain$ is given as candidates of the true distribution instead of the true distribution itself, we can extend the previous optimization problem to the problem introduced by Shtarkov \cite{Shtarkov:1987}:

\begin{equation}
\EqnLabel{ShtarkovMinMax}
\begin{split}
& \min_{\CodeLength{}{}} \max_{\DataSequence} \Bracket{\CodeLength{\DataSequence}{} - \min_{\DensityFunc{}{} \in \DensityDomain} \Paren{- \log \DensityFunc{\DataSequence}{}}} \\
& \ST \int \DataDiff \Exponential{- \CodeLength{\DataSequence}{}} \le 1.
\end{split}
\end{equation}

This problem is no longer trivial, and Shtarkov showed that the NML code length defined below attains its minimum \cite{Shtarkov:1987}:

\begin{equation}
\begin{split}
\NMLCode{\DataSequence} & \DefEq - \log \NMLDensityFunc{\DataSequence} \\
& = - \log \max_{\DensityFunc{}{} \in \DensityDomain} \DensityFunc{\DataSequence}{} + \log \int_{\DataDomain^\DataCnt} \max_{\DensityFunc{}{} \in \DensityDomain} \DensityFunc{\DataSequence}{} \Diff \DataSequence.
\end{split}
\end{equation}

The problem \EqnRef{ShtarkovMinMax} is reduced to \EqnRef{ShannonMinMax} if $\DensityDomain = \Set{\TrueDensityFunc{}}$, and in this sense, the NML code length is an extension of the self-entropy.

The NML code is one of the universally optimal codings when the true distribution in the given set is unknown \cite{Rissanen:1989}. 

The NML code length is widely used in model selection on the basis of the minimum description length principle \cite{Rissanen+:1992} \cite{Yamanishi:1992} \cite{Rissanen:1995} \cite{Rissanen:1997}. 
Here, the model that minimizes the NML code length for given data is selected.
Recently, it is shown that the NML code length bounds the generalized loss \cite{Grunwald+:2017}.

The calculation of the NML code length has been an important problem. Rissanen derived an asymptotic formula of the NML code length \cite{Rissanen:1996}, which clarified the behavior of the NML code length with $o(1)$ terms excluded as follows:
\begin{equation}
\begin{split}
& \log \int \DataDiff \DensityFunc{\DataSequence}{\MLE{\ParameterVec}{\DataSequence}} \\
& =
\frac{1}{2} \ParameterCnt \log \frac{\DataCnt}{2 \PiUnit}  + \log \int_{\ProperParameterDomain} \ParameterDiff \sqrt{\Det{\FisherMat{\ParameterVec}}} + \SmallOrder{1},
\end{split}
\end{equation}
where $\ParameterCnt$ denotes the dimension of the parameter.
This formula holds with certain regularity conditions and does not depend on the details of the model.

According to this formula, we can apply Nishii's analysis in terms of consistency in the selected model \cite{Nishii:1988} and Barron and Cover's result in terms of statistical risk \cite{Barron+:1991} to the model selection using the NML code length.

In contrast to the generality of Rissanen's asymptotic formula, a non-asymptotic calculation formula has been derived through model-by-model discussion \cite{Rissanen:2000} \cite{Kontkanen+:2007a} \cite{Kontkanen+:2007b} \cite{Kontkanen+:2008} \cite{Roos+:2008}.
Recently, Hirai and Yamanishi non-asymptotically calculated the NML code length for several models in the exponential family \cite{Hirai+:2013}.
They reduced the calculation of the NML code length to the parameter domain integral of the function denoted by $g$ in the paper.
However, the method to obtain Hirai-Yamanishi's $g$-function explicitly depends on the model. Thus, the exact calculation of the NML code length has been limited to particular models.

This paper proposes a novel methodology for calculating the NML code length on the basis of Fourier transformation. Our methodology enables the systematic analysis of NML code length in terms of both asymptotic expansion and exact calculation. As corollaries, our methodology provides an asymptotic formula with weaker assumptions compared to Rissanen's and a useful exact calculation formula for the exponential family.


\subsection{Significance of This Paper}
This paper proposes an alternative form of the NML code length based on the Fourier transform. Our form enables the systematic calculation of the NML code length. Specifically, it results in the two formulae presented below.
\subsubsection{Asymptotic Formula with Weaker Assumption}
Taking the limitation of our form leads to Rissanen's asymptotic formula \cite{Rissanen:1996}.
It should be noted that Lebesgue's dominant convergence theorem can be applied to our Fourier-transform-based form, which results in an asymptotic formula with a weaker assumption compared to that in the original paper \cite{Rissanen:1996}.

\subsubsection{Exact Calculation Formula for Exponential Family}
Our Fourier-transform-based form gives a simple formula for the exact calculation of the NML code length of the exponential family. The formula yields the NML code length from the partition function and the relationship between the canonical parameters and the expectation of sufficient statistics.


\subsection{Related Work}
\subsubsection{Asymptotic Formula with Weaker Assumption}
The consequence of the asymptotic formula in this paper is the same as that of Rissanen's theorem \cite{Rissanen:1996}. However, Rissanen's theorem assumes both the uniform asymptotic normality of the maximum likelihood estimator as well as the existence of a non-zero lower bound and an finite upper bound of the Fisher information; in contrast, our theorem does not involve these assumptions and allows the Fisher information to converge to zero or diverge towards the boundary.

\subsubsection{Exact Calculation Formula for Exponential Family}
Hirai and Yamanishi presented the exact calculation of several models in the exponential family through the integral of the $g$-function. However, in general, it is still difficult to obtain the explicit form of the $g$-function. In this paper, the general exact calculation formula for the exponential family is obtained, including Hirai and Yamanishi's results.

\section{Normalized Maximum Likelihood Code Length}
We consider a sequence $\DataSequence \DefEq \Range{\DataVec_{1}}{\DataVec_{2}}{\DataVec_{\DataCnt}}$ of continuous random variables and assume that they have a probability density function. 
\begin{definition}
Let $\DensityDomain \subset \SetBuilder{\DensityFunc{}{}: \Real^{\DataDimCnt \times \DataCnt} \supset \DataDomain^\DataCnt \to [0, \infty)}{\int_{\DataDomain^\DataCnt} \DensityFunc{\DataSequence}{} \Diff \DataSequence = 1}$ denote a set of density functions. Here, $\DataDomain \subset \Real^\DataDimCnt$ denotes the domain of a datum. Assume that $\max_{\DensityFunc{}{} \in \DensityDomain} \DensityFunc{\DataSequence}{}$ is a measurable function of $\DataSequence$. 
The NML code length is defined as its negative log likelihood as follows:
\begin{equation}
\begin{split}
\NMLCode{\DataSequence}
& \DefEq - \log \max_{\DensityFunc{}{} \in \DensityDomain} \DensityFunc{\DataSequence}{} + \log \MLComplexity{\DensityDomain},
\end{split}
\end{equation}
where $\MLComplexity{\DensityDomain} \DefEq \int_{\DataDomain^\DataCnt} \DataDiff \max_{\DensityFunc{}{} \in \DensityDomain} \DensityFunc{\DataSequence}{}$.
\end{definition}
$\MLComplexity{\DensityDomain}$ (or its logarithm) is called the \emph{parametric complexity (PC)} of $\DensityDomain$.
In this paper, we focus on the case in which it is easy to evaluate the first term but difficult to evaluate the maximum likelihood complexity. This is because, when even the first term is intractable, it is hardly possible to strictly evaluate the second term. 
In this paper, we consider the independent identical parametric model:
\begin{equation}
\DensityDomain_{\ParameterDomain} = \SetBuilder{\DensityFunc{}{}: \DataDomain^\DataCnt \to [0, \infty)}{\begin{array}{l}\DensityFunc{\DataSequence}{} = \DataProd \DensityFunc{\DataVec_\DataIdx}{\ParameterVec}, \\ \ParameterVec \in \ParameterDomain \subset \Real^{\ParameterCnt}\end{array}}
\end{equation}
as a set of density functions. Here, $\ParameterVec$ is its parameter and $\ParameterDomain$ is the domain of the parameter.

We mainly analyze a proper parameter domain $\ProperParameterDomain \subset \ParameterDomain$ defined as follows.
\begin{definition}
A subset $\ProperParameterDomain$ of $\ParameterDomain$ is proper if the following conditions are satisfied:
\begin{enumerate}
\item Map $\ProperParameterDomain \ni \ParameterVec \mapsto \DensityFunc{\cdot}{\ParameterVec} \in \ProperDensityDomain \subset \DensityDomain$ is bijective (one to one).
\item For all $\DataSequence \in \DataDomain^\DataCnt$, a unique solution $\MLE{\ParameterVec}{\DataSequence}$ of $\max_{\ParameterVec \in \ProperParameterDomain} \DensityFunc{\DataSequence}{\ParameterVec}$ exists; that is, a unique maximum likelihood estimator (MLE) $\MLE{\ParameterVec}{\DataSequence}$ exists.
\item $\max_{\ParameterVec \in \ProperParameterDomain} \DensityFunc{\DataSequence}{\ParameterVec}$ is a measurable function of $\DataSequence$.
\item If $\ParameterVec \in \ProperParameterDomain$ and $\DataVec_{\DataIdx} \DistributedAs \DensityFunc{\DataVec_\DataIdx}{\ParameterVec}$, the asymptotic normality of the MLE $\MLE{\ParameterVec}{\DataSequence}$ holds; that is, $\sqrt{\DataCnt} \Paren{\MLE{\ParameterVec}{\DataSequence} - \ParameterVec} \WeakConvergence \NormalDist{\ZeroVec}{\FisherMat{\ParameterVec}^\Inverse}$, where $\FisherMat{\ParameterVec}^\Inverse$ denotes the Fisher information matrix.
\end{enumerate}
We also define the proper data sequence domain $\ProperDataSequenceDomain$ as follows:
\begin{equation}
\ProperDataSequenceDomain \DefEq \SetBuilder{\DataSequence \in \DataSequenceDomain}{\max_{\ParameterVec \in \ProperParameterDomain} \DensityFunc{\DataSequence}{\ParameterVec} = \max_{\ParameterVec \in \ParameterDomain} \DensityFunc{\DataSequence}{\ParameterVec}}.
\end{equation}
\end{definition}
\begin{remark}
Sufficient conditions for 4) have been discussed (for example, see \cite{Van:1998}).
At least the positive definiteness of $\FisherMat{\ParameterVec}^\Inverse$ in $\ProperDataDomain$ is required for 4).
\end{remark}
\begin{remark}
Since $\ProperParameterDomain \subset \ParameterDomain$, the following holds in general: $\max_{\ParameterVec \in \ProperParameterDomain} \DensityFunc{\DataSequence}{\ParameterVec} \le \max_{\ParameterVec \in \ParameterDomain} \DensityFunc{\DataSequence}{\ParameterVec}$.
\end{remark}
\begin{remark}
In this paper, $\MLE{\ParameterVec}{\DataSequence}$ always denotes the unique MLE on $\ProperParameterDomain$. If $\ProperParameterDomain \subsetneq \ParameterDomain$, the MLE in $\ParameterDomain$ can be non-unique. 
\end{remark}
Roughly speaking, the proper parameter domain is a tractable subset of the model, and the proper data sequence domain is a set of sequences, the MLE of which lies in the proper parameter domain. 
The PC can be decomposed as follows:
\begin{equation}
\begin{split}
\EqnLabel{Decomposition}
& \int_{\DataSequenceDomain} \DataDiff \max_{\ParameterVec \in \ParameterDomain} \DensityFunc{\DataSequence}{\ParameterVec} = \int_{\DataSequenceDomain} \DataDiff \DensityFunc{\DataSequence}{\MLE{\ParameterVec}{\DataSequence}} \\
& + \int_{\DataSequenceDomain \setminus \ProperDataSequenceDomain} \DataDiff \Paren{\max_{\ParameterVec \in \ParameterDomain} \DensityFunc{\DataSequence}{\ParameterVec} - \DensityFunc{\DataSequence}{\MLE{\ParameterVec}{\DataSequence}}}. \\
\end{split}
\end{equation}
\begin{remark}
If we can take $\ProperParameterDomain = \ParameterDomain$ (as is often the case with a well-behaved model such as the exponential family models), the second term vanishes, and the logarithm of the first term is equivalent to the parametric complexity.
\end{remark}
We assume that the second term is ignorable and focus on the first term $\MLComplexity{\ProperParameterDomain} \DefEq \int_{\DataSequenceDomain} \DataDiff \DensityFunc{\DataSequence}{\MLE{\ParameterVec}{\DataSequence}}$ in this paper. 

Note that $\MLComplexity{\ProperParameterDomain}$ carries excessive data sequences and often diverges to infinity. To avoid this problem, we introduce \emph{luckiness} \cite{Grunwald:2007} to generalize $\MLComplexity{\ProperParameterDomain}$ as follows:

\begin{definition}
Let $\ParameterWeightFunc{}: \ParameterDomain \to [0 , \infty)$ denote a weight function on $\ProperParameterDomain$ called \emph{luckiness}.
We define the \emph{luckiness parametric complexity (LPC)} of $\ProperParameterDomain$ as follows:
\begin{equation}
\begin{split}
\LMLComplexity{\ParameterWeightFunc{}}{\ProperParameterDomain} \DefEq \int \DataDiff \DensityFunc{\DataSequence}{\MLE{\ParameterVec}{\DataSequence}} \ParameterWeightFunc{\MLE{\ParameterVec}{\DataSequence}},
\end{split}
\end{equation}
where $\MLE{\ParameterVec}{\DataSequence} \DefEq \ArgMax_{\ParameterVec \in \ProperParameterDomain} \DensityFunc{\DataSequence}{\ParameterVec}$.
\end{definition}

\begin{remark}
If $\ParameterWeightFunc{\ParameterVec} = 1$, the LPC is equivalent to the PC.
\end{remark}
Let $A$ be a subset of $\ProperParameterDomain$. We can regard the LPC $\LMLComplexity{\Indicator{A}}{\ProperParameterDomain}$ as a restriction of the PC $\MLComplexity{\ProperParameterDomain}$  to $A$, where $\Indicator{\cdot}$ denotes the indicator function. This restriction is often necessary and used in continuous variable cases \cite{Grunwald:2007} \cite{Hirai+:2013}.

\section{Fourier Form of NML Code Length}
First, we make assumptions that allow us to exchange integrals.

\begin{assumption}
\AssumptionLabel{IntIntExchange}
\begin{enumerate}
\item For all $\varPhi_0 \subset \ProperParameterDomain$ that have measure zero, $\SetBuilder{\DataSequence}{\MLE{\ParameterVec}{\DataSequence} \in \varPhi_0}$ has measure zero.
\item For all $\DataSequence$, $\DensityFunc{\DataSequence}{\ParameterVec} \ParameterWeightFunc{\ParameterVec}$
is integrable and square-integrable
as a function of $\ParameterVec$.
\item For all $\DataSequence$, the Fourier transform $\DensityWeightFT{\DataSequence}{\FrequencyVec}$
of  $\DensityFunc{\DataSequence}{\ParameterVec} \ParameterWeightFunc{\ParameterVec}$ is integrable as a function of $\FrequencyVec$,
where
\begin{equation}
\begin{split}
&\DensityWeightFT{\DataSequence}{\FrequencyVec} \\ 
&\DefEq \Paren{\frac{1}{2 \PiUnit}}^{\frac{\ParameterCnt}{2}} \int_{\ProperParameterDomain} \ParameterDiff \Exponential{- \IUnit \FrequencyVec^{\Transpose} \ParameterVec} \DensityFunc{\DataSequence}{\ParameterVec} \ParameterWeightFunc{\ParameterVec}.
\end{split}
\end{equation}
\item The characteristic function $\MLECharacteristic{\ParameterVec}{\DataCnt}{\FrequencyVec}$ of the maximum likelihood estimator is integrable as a function of ${\ParameterVec}$ and ${\FrequencyVec}$,
where
\begin{equation}
\MLECharacteristic{\ParameterVec}{\DataCnt}{\FrequencyVec} \DefEq \int \DataDiff \DensityFunc{\DataSequence}{\ParameterVec} \Exponential{\IUnit \FrequencyVec^{\Transpose} \Paren{\MLE{\ParameterVec}{\DataSequence} - \ParameterVec}}. 
\end{equation}
\end{enumerate}
\end{assumption}
We obtain the Fourier-transform-based form of the NML code length as follows: 
\begin{theorem}
\ThmLabel{Basic}
Under \AssumptionRef{IntIntExchange}, the PC is calculated as follows:
\begin{equation}
\begin{split}
& \int \DataDiff \DensityFunc{\DataSequence}{\MLE{\ParameterVec}{\DataSequence}} \ParameterWeightFunc{\MLE{\ParameterVec}{\DataSequence}}
= \int \ParameterDiff \ParameterWeightFunc{\ParameterVec} \ASGFunc{\DataCnt}{\ParameterVec}, 
%
\end{split}
\end{equation}
where
\begin{equation}
\begin{split}
& \ASGFunc{\DataCnt}{\ParameterVec} \\
& \DefEq \frac{1}{\Paren{2 \PiUnit}^{\ParameterCnt}} \int \FrequencyDiff \int \DataDiff \DensityFunc{\DataSequence}{\ParameterVec} \Exponential{\IUnit \FrequencyVec^{\Transpose} \Paren{\MLE{\ParameterVec}{\DataSequence} - \ParameterVec}}
\end{split}
\end{equation}
\end{theorem}



\begin{IEEEproof}
Since $\DensityFunc{\DataSequence}{\ParameterVec} \in \LP{1}{\CompactParameterDomain} \cap \LP{2}{\CompactParameterDomain}$, 
\begin{equation}
\begin{split}
& \DensityFunc{\DataSequence}{\ParameterVec'} \ParameterWeightFunc{\MLE{\ParameterVec}{\DataSequence}} \\
& = \Paren{\frac{1}{2 \PiUnit}}^{\frac{\ParameterCnt}{2}} \int \FrequencyDiff \Exponential{\IUnit \FrequencyVec^{\Transpose} \ParameterVec'} \DensityWeightFT{\DataSequence}{\FrequencyVec} \quad \AS.
\end{split}
\end{equation}
Thus, the following holds with \AssumptionRef{IntIntExchange} 1):
\begin{equation}
\begin{split}
&\int \DataDiff \DensityFunc{\DataSequence}{\MLE{\ParameterVec}{\DataSequence}} \ParameterWeightFunc{\MLE{\ParameterVec}{\DataSequence}} \\
&= \int \DataDiff
\Paren{\frac{1}{2 \PiUnit}}^{\frac{\ParameterCnt}{2}} \int \FrequencyDiff \Exponential{\IUnit \FrequencyVec^{\Transpose} \MLE{\ParameterVec}{\DataSequence}} \DensityWeightFT{\DataSequence}{\FrequencyVec} \\
&= \Paren{\frac{1}{2 \PiUnit}}^{\frac{\ParameterCnt}{2}} \int \FrequencyDiff \int \DataDiff \Exponential{\IUnit \FrequencyVec^{\Transpose} \MLE{\ParameterVec}{\DataSequence}} \DensityWeightFT{\DataSequence}{\FrequencyVec},
\end{split}
\end{equation}
where the last equation follows from the (absolute) integrability of $\DensityWeightFT{\DataSequence}{\FrequencyVec}$ and Fubini's theorem.
We exchange integrals again likewise as follows:
\begin{equation}
\begin{split}
& \int_{\ProperDataSequenceDomain} \DataDiff \Exponential{\IUnit \FrequencyVec^{\Transpose} \MLE{\ParameterVec}{\DataSequence}} \DensityWeightFT{\DataSequence}{\FrequencyVec} \\
&= \Paren{\frac{1}{2 \PiUnit}}^{\frac{\ParameterCnt}{2}} \int \DataDiff \Exponential{\IUnit \FrequencyVec^{\Transpose} \MLE{\ParameterVec}{\DataSequence}} \\
& \quad \times \int \ParameterDiff \Exponential{- \IUnit \FrequencyVec^{\Transpose} \ParameterVec} \DensityFunc{\DataSequence}{\ParameterVec} \ParameterWeightFunc{\ParameterVec} \\
&= \Paren{\frac{1}{2 \PiUnit}}^{\frac{\ParameterCnt}{2}} \int \ParameterDiff \ParameterWeightFunc{\ParameterVec} \\
& \quad \times \int \DataDiff \Exponential{\IUnit \FrequencyVec^{\Transpose} \Paren{\MLE{\ParameterVec}{\DataSequence} - \ParameterVec}} \DensityFunc{\DataSequence}{\ParameterVec},
\end{split}
\end{equation}
where the last equation follows from the (absolute) integrability of $\DensityWeightFT{\DataSequence}{\FrequencyVec}$ and Fubini's theorem.
By the third assumption, we can exchange the integral with respect to $\ParameterVec$ and $\FrequencyVec$, which completes the proof.
\end{IEEEproof}

\section{Asymptotic Formula}
By taking the limitation of \ThmRef{Basic}, we can prove the asymptotic formula of the LPC, which relaxes some conditions given by Rissanen's asymptotic formula \cite{Rissanen:1996}.
First, we make assumptions that allow us to exchange the limitation and integral.
\begin{assumption}
\AssumptionLabel{LimIntExchange}
\begin{enumerate}
\item There exists an integrable function $\MLECharacteristicBound{\ParameterVec}{\FrequencyVec}$ of $\FrequencyVec$ such that $\Abs{\MLECharacteristic{\ParameterVec}{\DataCnt}{\FrequencyVec}} < \MLECharacteristicBound{\ParameterVec}{\FrequencyVec}$, for all $\DataCnt$ and $\ParameterVec$.
\item There exists an integrable function $\ASGFuncBound{\ParameterVec}$ of $\ParameterVec$ such that $\Abs{\ASGFunc{\DataCnt}{\ParameterVec}} < \ASGFuncBound{\ParameterVec}$.
\end{enumerate}
\end{assumption}
\begin{theorem}[Asymptotic formula of the NML code length]
Under \AssumptionRef{IntIntExchange} and \AssumptionRef{LimIntExchange},
the following holds:
\begin{equation}
\begin{split}
&\log \int \DataDiff \DensityFunc{\DataSequence}{\MLE{\ParameterVec}{\DataSequence}} \ParameterWeightFunc{\MLE{\ParameterVec}{\DataSequence}} \\
& =
\frac{1}{2} \ParameterCnt \log \frac{\DataCnt}{2 \PiUnit} + \log \int_{\ProperParameterDomain} \ParameterDiff \ParameterWeightFunc{\ParameterVec} \sqrt{\Det{\FisherMat{\ParameterVec}}}
+ \SmallOrder{1}.
\end{split}
\end{equation}
\end{theorem}



\begin{IEEEproof}
The assumptions allow us to apply Lebesgue's dominant convergence theorem to \ThmRef{Basic} as follows: 
\begin{equation}
\begin{split}
& \lim_{\DataCnt \to \infty} \int \DataDiff \DensityFunc{\DataSequence}{\MLE{\ParameterVec}{\DataSequence}} \ParameterWeightFunc{\MLE{\ParameterVec}{\DataSequence}} \\
& = \int \ParameterDiff \ParameterWeightFunc{\ParameterVec} \lim_{\DataCnt \to \infty} \ASGFunc{\DataCnt}{\ParameterVec} \\
& = \int \ParameterDiff \ParameterWeightFunc{\ParameterVec} \frac{1}{\Paren{2 \PiUnit}^{\ParameterCnt}} \int \FrequencyDiff \lim_{\DataCnt \to \infty} \MLECharacteristic{\ParameterVec}{\DataCnt}{\FrequencyVec} 
\end{split}
\end{equation}
If $\Parameter \in \Interior{\ProperParameterDomain}$,
\begin{equation}
\sqrt{\DataCnt} \Paren{\MLE{\ParameterVec}{\DataSequence} - \ParameterVec}
\WeakConvergence
\NormalDist{\ZeroVec}{\FisherMat{\ParameterVec}^\Inverse},
\end{equation}
by asymptotic normality of the MLE.
Hence, by Levy's continuity theorem,
\begin{equation}
\begin{split}
\lim_{\DataCnt \to \infty} \MLECharacteristic{\ParameterVec}{\DataCnt}{\FrequencyVec} = \Exponential{- \frac{1}{2} \FrequencyVec^{\Transpose} \FisherMat{\ParameterVec} \FrequencyVec}^{\Inverse},
\end{split}
\end{equation}
which completes the proof.
\end{IEEEproof}


\begin{remark}
The consequence of the theorem is the same as that of Rissanen's formula \cite{Rissanen:1996}. In contrast to Rissanen's formula, we do not make assumptions on the boundedness of the determinant of the Fisher information matrix or the uniform asymptotic normality of the MLE. Thus, we can expect that our formula is easy to apply even when $\ProperParameterDomain$ is not compact and the boundedness of the Fisher information and uniform asymptotic normality of the MLE are difficult to guarantee. 
\end{remark}


\begin{table*}[hbtp]
  \TabLabel{Exponential}
  \caption{Exponential Family and PC}
  \centering
  \begin{tabular}{lllllll}
    \hline
Distribution & Density
& Sufficient & Canonical parameter $\NaturalParameterVec$ & Partition & Parametric
\\
 & function
& statistics & & function & complexity
\\
 & $\DensityFunc{\DataSequence}{\ParameterVec}$
& $\StatisticElement{\ParameterIdx}{\DataVec}$ & Expectation parameter $\ExpectationParameterVec = \Expect{\StatisticElement{\ParameterIdx}{\DataVec}}$ & $\Partition{\NaturalParameterVec}$ \\
\hline 
\hline
Normal dist. with &     $\DensityFunc{\DataSequence}{\ExponentialMean}$ & $\Data$ & $\NaturalParameter = \frac{\NormalMean}{\NormalVariance} \in (- \infty, + \infty)$  & $\NormalVariance^{\frac{1}{2}} \Exponential{\frac{1}{2} \NormalVariance \NaturalParameter^{2}}$ & $\int_{-\infty}^{+\infty} \Diff \NormalMean \frac{\ParameterWeightFunc{\NormalMean}}{\sqrt{2 \PiUnit \frac{\NormalVariance}{\DataCnt}}}$
\\
\quad known variance $\NormalVariance$ & $= \frac{1}{\sqrt{2 \PiUnit \NormalVariance}} \Exponential{- \frac{\Paren{\Data - \NormalMean}^{2}}{2 \NormalVariance}}$ & & $ \NormalMean = \NormalVariance \NaturalParameter \in (- \infty, + \infty) $ & 
\\
\hline
Normal dist. with &     $\DensityFunc{\DataSequence}{\NormalVariance}$ & $\Paren{\Data - \mu}^2$ & $\NaturalParameter = - \frac{1}{2 \NormalVariance} \in (- \infty, 0) $ & $\frac{1}{\sqrt{- \NaturalParameter}}$ & $\frac{\Paren{\frac{1}{2} \DataCnt}^{\frac{1}{2} \DataCnt} \Exponential{- \frac{1}{2} \DataCnt}}{\GammaFunc{\frac{1}{2} \DataCnt}}$
\\
\quad known mean $\Mean$ & $= \frac{1}{\sqrt{2 \NormalVariance}} \Exponential{- \frac{\Paren{\Data - \Mean}^2}{2 \NormalVariance}}$ &  & $ \NormalVariance = - \frac{1}{2 \NaturalParameter} \in (0, \infty) $ & & $\times \int_{0}^{+\infty} \Diff \NormalVariance \frac{\ParameterWeightFunc{\NormalVariance}}{\NormalVariance}$
\\
\\
Laplace dist. with &     $\DensityFunc{\DataSequence}{\LaplaceScale}$ & $\Abs{\Data - \mu}$ & $\NaturalParameter = - \frac{1}{\LaplaceScale} \in (- \infty, 0) $ & $\frac{2}{- \NaturalParameter}$ & $\frac{\DataCnt^{\DataCnt} \Exponential{- \DataCnt}}{\GammaFunc{\DataCnt}}$
\\
\quad known mean $\Mean$ & $= \frac{1}{2 \LaplaceScale} \Exponential{- \frac{\Abs{\Data - \Mean}}{\LaplaceScale}}$ &  & $ \LaplaceScale = - \frac{1}{\NaturalParameter} \in (0, \infty) $ & & $\times \int_{0}^{+\infty} \Diff \LaplaceScale \frac{\ParameterWeightFunc{\LaplaceScale}}{\LaplaceScale}$
\\
\\
Gamma dist. with &     $\DensityFunc{\DataSequence}{\GammaMean}$ & $\Data$ & $\NaturalParameter = - \frac{\GammaShape}{\GammaMean} \in (- \infty, 0) $ & $\frac{1}{\Paren{- \NaturalParameter}^{\GammaShape}}$ & $\frac{\Paren{\GammaShape \DataCnt}^{\GammaShape \DataCnt} \Exponential{- \GammaShape \DataCnt}}{\GammaFunc{\GammaShape \DataCnt}}$
\\
\quad known shape $\GammaShape$ \footnotemark[1] \footnotemark[2] & $= \frac{\GammaShape^{\GammaShape} \Data^{\GammaShape-1}}{\GammaFunc{\GammaShape} \GammaMean^{\GammaShape}} \Exponential{- \frac{\GammaShape \Data}{\GammaMean}}$ & & $ \GammaMean = - \frac{\GammaShape}{\NaturalParameter} \in (0, \infty) $ & & $\times \int_{0}^{+\infty} \Diff \GammaMean \frac{\ParameterWeightFunc{\GammaMean}}{\GammaMean}$
\\
\\
Weibull dist. with &     $\DensityFunc{\DataSequence}{\WeibullExpectation}$ & $\Data^\WeibullShape$ & $\NaturalParameter = - \frac{1}{\WeibullExpectation} \in (- \infty, 0) $ &  $\frac{1}{- \NaturalParameter}$ & $\frac{\DataCnt^{\DataCnt} \Exponential{- \DataCnt}}{\GammaFunc{\DataCnt}}$
\\
\quad known shape $\WeibullShape$ & $= \WeibullShape \WeibullExpectation^{-\frac{\WeibullShape+1}{\WeibullShape}} \Data^{\WeibullShape} \Exponential{- \frac{\Data^{\WeibullShape}}{\WeibullExpectation}}$ & & $ \WeibullExpectation = - \frac{1}{\NaturalParameter} \in (0, \infty) $ & & $\times \int_{0}^{+\infty} \Diff \WeibullExpectation \frac{\ParameterWeightFunc{\WeibullExpectation}}{\WeibullExpectation}$
\\
\hline
Gamma dist. with &     $\DensityFunc{\DataSequence}{\NaturalParameter}$ & $\log \Data$ & $\NaturalParameter = \InverseDigammaFunc{\GammaLogExpectation - \log \GammaScale} - 1 \in (- 1, + \infty) $ & $\GammaFunc{\NaturalParameter+1} \GammaScale^{\NaturalParameter+1}$ & 
See \EqnRef{GammaShape}
\\
\quad known scale $\GammaScale$ \footnotemark[3] & $= \frac{\Data^{\NaturalParameter}}{\GammaFunc{\NaturalParameter+1} \GammaScale^{\NaturalParameter+1}} \Exponential{- \frac{\Data}{\GammaScale}}$ & & $ \GammaLogExpectation = - \DigammaFunc{\NaturalParameter+1} + \log \GammaScale \in (-\infty, +\infty) $ & 
\\
\hline
  \end{tabular}
\end{table*}

\section{Non-asymptotic Formula for Exponential Family}
First, we present the notation of the exponential family. Then, we present the non-asymptotic formula of the PC for the exponential family.
\subsection{Exponential Family and Its Canonical Parameters and Expectation Parameters}
We say that a model is in the exponential family when we can express its density function with its canonical parameters $\NaturalParameterVec \in \NaturalParameterDomain \subset \Real^{\ParameterCnt}$, sufficient statistics $\StatisticVec{}: \Real^{\DataDimCnt} \supset \DataDomain \to \Real^{\ParameterCnt}$, and base measure $\BaseMeasure{}: \Real^{\DataDimCnt} \supset \DataDomain \to [0, +\infty)$ as follows:
\begin{equation}
\DensityFunc{\DataVec}{\NaturalParameterVec} = \frac{\BaseMeasure{\DataVec}}{\Partition{\NaturalParameterVec}} \Exponential{\NaturalParameterVec^{\Transpose} \StatisticVec{\DataVec}}.
\end{equation}
Here, the partition function $\Partition{}: \NaturalParameterDomain \to [0, \infty)$ is defined by
$\Partition{\NaturalParameterVec} = \int \Diff \DataVec \BaseMeasure{\DataVec} \Exponential{\NaturalParameterVec^{\Transpose} \StatisticVec{\DataVec}}$.




We define the transform from the canonical parameters to expectation parameters as
\begin{equation}
\begin{split}
\ExpectationParameterElementMap{\ParameterIdx}{\NaturalParameterVec}
&\DefEq
\int \Diff \DataVec \StatisticElement{\ParameterIdx}{\DataVec} \frac{\BaseMeasure{\DataVec}}{\Partition{\NaturalParameterVec}} \Exponential{\ParameterSum \NaturalParameter_{\ParameterIdx} \StatisticElement{\ParameterIdx}{\DataVec}} 
, \\
\ExpectationParameterVecMap{\NaturalParameterVec}
&\DefEq
\BRowMat{\ExpectationParameterElementMap{1}{\NaturalParameterVec}}{\ExpectationParameterElementMap{2}{\NaturalParameterVec}}{\ExpectationParameterElementMap{\ParameterCnt}{\NaturalParameterVec}}^\Transpose
\end{split}
\end{equation}
and let $\NaturalParameterVecMap{\ExpectationParameterVec}$ denote its inverse transform, 
assuming $\ExpectationParameterVecMap{\cdot}$ is bijective.



We define the MLE with respect to the expectation parameters as follows:
\begin{equation}
\MLE{\ExpectationParameterVec}{\DataSequence}
\DefEq \ArgMax_{\ExpectationParameterVec} \DataProd \DensityFunc{\DataVec_\DataIdx}{\NaturalParameterVecMap{\ExpectationParameterVec}}.
\end{equation}

We can calculate the PC of a model in the exponential family as follows.

\begin{theorem}
\ThmLabel{Exponential}
Let $\DensityFunc{\DataSequence}{\NaturalParameterVec}$ be the density function of a model in the exponential family, where $\NaturalParameterVec$ denotes its natural parameter and $\ExpectationParameter$ denotes its expectation parameter.
The LPC of $\DensityFunc{\DataSequence}{\NaturalParameterVec}$ is expressed as follows:
\begin{equation}
\begin{split}
& \int \DataDiff \DensityFunc{\DataSequence}{\NaturalParameterVecMap{\MLE{\ExpectationParameterVec}{\DataSequence}}} \ParameterWeightFunc{\MLE{\ExpectationParameterVec}{\DataSequence}} \\
& = \Paren{\frac{1}{2 \PiUnit}}^{\ParameterCnt} \int \ExpectationParameterDiff \ParameterWeightFunc{\ExpectationParameterVec} \int \FrequencyDiff \Exponential{- \IUnit \FrequencyVec^\Transpose \ExpectationParameterVec} \Paren{\frac{\Partition{\NaturalParameterVecMap{\ExpectationParameterVec} + \IUnit \frac{\FrequencyVec}{\DataCnt}}}{\Partition{\NaturalParameterVecMap{\ExpectationParameterVec}}}}^\DataCnt 
.
\end{split}
\end{equation}
\end{theorem}

\begin{corollary}
\CorLabel{Exponential}
Let $\Range{X_{1}^{\Paren{\DataCnt, \ExpectationParameterVec}}}{X_{2}^{\Paren{\DataCnt, \ExpectationParameterVec}}}{X_{\DataCnt}^{\Paren{\DataCnt, \ExpectationParameterVec}}}$ be a sequence of i.i.d. $\ParameterCnt$-dimensional random variables, the characteristic function of each of which is given by 
$\frac{\Partition{\NaturalParameterVecMap{\ExpectationParameterVec} + \IUnit \frac{\TmpFrequencyVec}{\DataCnt}}}{\Partition{\NaturalParameterVecMap{\ExpectationParameterVec}}}$,
and let $g^{\Paren{\DataCnt}} \Paren{\DataVec, \ExpectationParameterVec}$ be the density function of $\DataSum X_{\DataIdx}^{\Paren{\DataCnt, \ExpectationParameterVec}}$.
Then, the PC is expressed as follows:
\begin{equation}
\begin{split}
& \int \DataDiff \DensityFunc{\DataSequence}{\NaturalParameterVecMap{\MLE{\ExpectationParameterVec}{\DataSequence}}} \ParameterWeightFunc{\MLE{\ExpectationParameterVec}{\DataSequence}} \\
& = \Paren{\frac{1}{2 \PiUnit}}^{\ParameterCnt} \int \ExpectationParameterDiff g^{\Paren{\DataCnt}} \Paren{\ExpectationParameterVec, \ExpectationParameterVec} \ParameterWeightFunc{\ExpectationParameterVec}. \\
\end{split}
\end{equation}
\end{corollary}

\begin{remark}
\ThmRef{Exponential} reduces the original $\DataCnt$-times integrals to a $2 \ParameterCnt$-times integral. \CorRef{Exponential} implies that, if we know the density function of $\DataSum X_{\DataIdx}^{\Paren{\DataCnt, \ExpectationParameterVec}}$, the characteristic function of each of which is given using the partition function of the original function, the calculation of the original PC can be reduced to one integral calculation, and it is often analytically obtained. 
\end{remark}



\begin{IEEEproof}
Note that the following holds with respect to its maximum likelihood estimator as follows:
\begin{equation}
\begin{split}
\frac{\partial}{\partial \NaturalParameterVec} \log \Partition{\MLE{\NaturalParameterVec}{\DataSequence}}
& =
\frac{1}{\DataCnt} \DataSum \StatisticVec{\DataVec_\DataIdx}.
\end{split}
\end{equation}

Also note that the maximum likelihood estimator $\MLE{\ExpectationParameterVec}{\DataSequence}$ with respect to the expectation parameters can be written as $\MLE{\ExpectationParameterVec}{\DataSequence} = \ExpectationParameterVec \Paren{\MLE{\NaturalParameterVec}{\DataSequence}}$. Here it holds that
\begin{equation}
\MLEElement{\ExpectationParameter}{\ParameterIdx}{\DataSequence} = \frac{1}{\DataCnt} \DataSum \StatisticElement{\ParameterIdx}{\DataVec_\DataIdx}
\end{equation}

We can calculate $g$ function as follows:
\begin{equation}
\begin{split}
& \Paren{2 \PiUnit}^{\ParameterCnt} \ASGFunc{\DataCnt}{\ParameterVec} \\
& = \int \FrequencyDiff \int \DataDiff \left[ \Bracket{\DataProd \frac{\BaseMeasure{\DataVec_\DataIdx}}{\Partition{\NaturalParameterVec \Paren{\ExpectationParameterVec}}} \Exponential{\NaturalParameterVecMap{\ExpectationParameterVec}^{\Transpose} \StatisticVec{\DataVec_\DataIdx}}} \right. \\ 
& \quad \times \left. \Exponential{\IUnit \FrequencyVec^{\Transpose} \Paren{\frac{1}{\DataCnt} \DataSum \StatisticVec{\DataVec_\DataIdx} - \ExpectationParameterVec}} \right]
\\
& = \int \FrequencyDiff \Exponential{- \IUnit \FrequencyVec^{\Transpose} \ExpectationParameterVec} \int \DataDiff \left[ \Paren{\frac{\BaseMeasure{\DataVec_\DataIdx}}{\Partition{\NaturalParameterVecMap{\ExpectationParameterVec}}}}^{\DataCnt} \right. \\ 
& \quad \times \left. \Exponential{\DataSum \Paren{\NaturalParameterVecMap{\ExpectationParameterVec} + \IUnit \frac{\FrequencyVec}{\DataCnt}}^{\Transpose} \StatisticVec{\DataVec_\DataIdx}} 
\right], 
\end{split}
\end{equation}
Substituting this to \ThmRef{Basic} completes the proof.
\end{IEEEproof}

\footnotetext[1]{Including the exponential distribution with $\GammaShape = 1$}
\footnotetext[2]{Derived in \cite{Hirai+:2013}}
\footnotetext[3]{Including the chi-squared distribution with $\GammaScale = 2$}

\subsection{Examples}
In this subsection, we give examples of PC (LPC) calculation using \ThmRef{Exponential}. These examples include the results in \cite{Hirai+:2013}. \TabRef{Exponential} lists the results.

\subsubsection{Fixed Variance Distribution}

If the relationship between the natural parameter and expectation parameter is given by $\NaturalParameter = \frac{\ExpectationParameter}{v}$ with a constant $v$ and the partition function is given by $\Partition{\NaturalParameter} = C \Exponential{\frac{v^{2} \NaturalParameter^{2}}{2 D}}$ with a constant $D$, we can calculate the LPC as follows: 

\begin{equation}
\begin{split}
& \Paren{\frac{1}{2 \PiUnit}} \int_{-\infty}^{+\infty} \Diff \ExpectationParameter \left[ \ParameterWeightFunc{\ExpectationParameter} \right. \\
& \quad \left. \times \int_{-\infty}^{+\infty} \Diff \Frequency \Exponential{\frac{\DataCnt}{D} \Paren{\frac{\ExpectationParameter^2}{2} - \frac{\Paren{\ExpectationParameter + \IUnit \frac{\Frequency}{\DataCnt}}^2}{2}}} \Exponential{- \IUnit \Frequency \ExpectationParameter} \right] \\
& = \int_{-\infty}^{+\infty} \Diff \ExpectationParameter \frac{\ParameterWeightFunc{\ExpectationParameter}}{\sqrt{2 \PiUnit \frac{D}{\DataCnt}}}. \\
\end{split}
\end{equation}




\subsubsection{Exponential Distribution Type}
If the relationship between the natural parameter and expectation parameter is given by $\NaturalParameter=-\frac{C}{\ExpectationParameter}$ with a constant $C, m$ and the partition function is given by $\Partition{\NaturalParameter} = \frac{D}{\Paren{- \NaturalParameter}^{m}}$ with a constant $D$, we can calculate the LPC as follows: 
\begin{equation}
\begin{split}
& \Paren{\frac{1}{2 \PiUnit}} \int_{0}^{+\infty} \Diff \ExpectationParameter \ParameterWeightFunc{\ExpectationParameter} \int \Diff \Frequency \Paren{\frac{- \frac{C}{\ExpectationParameter}}{- \frac{C}{\ExpectationParameter} + \IUnit \frac{\Frequency}{\DataCnt}}}^{m \DataCnt} \Exponential{ - \IUnit \Frequency \ExpectationParameter} \\
& = \frac{\Paren{C \DataCnt}^{m \DataCnt} \Exponential{- C \DataCnt}}{\GammaFunc{m \DataCnt}} \int_{0}^{+\infty} \Diff \ExpectationParameter \frac{\ParameterWeightFunc{\ExpectationParameter}}{\ExpectationParameter}.  
\end{split}
\end{equation}
The last equation holds because $\Paren{\frac{- \frac{C}{\ExpectationParameter}}{- \frac{C}{\ExpectationParameter} + \IUnit \frac{\Frequency}{\DataCnt}}}^{m \DataCnt}$ is the characteristic function of the gamma distribution with the shape parameter $m \DataCnt$ and rate parameter $\frac{\ExpectationParameter}{C \DataCnt}$.
If we set $\ParameterWeightFunc{\ExpectationParameter} = \Indicator{[\ExpectationParameter_{\min}, \ExpectationParameter_{\max}]} \Paren{\ExpectationParameter}$, the PC is equal to $\frac{\DataCnt^{m \DataCnt} \Exponential{- \DataCnt}}{\GammaFunc{m \DataCnt}} \log{\frac{\ExpectationParameter_{\max}}{\ExpectationParameter_{\min}}}$.
This formula can be applied to distributions including the exponential distribution, chi-squared distribution, Laplace distribution with a known mean, Weibull distribution with a known shape, and gamma distribution with a known shape.

\subsubsection{Chi-squared Distribution Type}
We discuss the gamma distribution with a known scale $\GammaScale$. The result here includes the chi-squared distribution (set $\GammaScale = 2$). we can calculate the LPC as follows: 
\begin{equation}
\EqnLabel{GammaShape}
\begin{split}
\\
& \int_{-\infty}^{+\infty} \Diff s \left[ \ParameterWeightFunc{\log \GammaScale - s} \int \Diff \TmpFrequency \Paren{\frac{\GammaFunc{\InverseDigammaFunc{-s} - \IUnit \frac{\TmpFrequency}{\DataCnt}}}{\GammaFunc{\InverseDigammaFunc{-s}}}}^{\DataCnt} \right. \\
& \quad \times \left. \Exponential{- \IUnit \TmpFrequency s} \right]
\\
& = \int_{-\infty}^{+\infty} \Diff s \ParameterWeightFunc{\log \GammaScale - s} G^{(\DataCnt)} \Paren{s; \InverseDigammaFunc{-s}, \frac{1}{\DataCnt}}.
\end{split}
\end{equation}
Here, $G^{(\DataCnt)} \Paren{\Data; p, q}$ is the probability density function of the sum of $\DataCnt$ i.i.d. samples, the density of which is given by
\begin{equation}
G \Paren{\Data; p, q} = \frac{1}{q \GammaFunc{p}}\Exponential{- \frac{p}{q} \Data - \Exponential{- \frac{\Data}{q}}}.
\end{equation}

\section{Conclusion}
In this paper, we derived a non-asymptotic form of the NML code length and clarified its relationship to the asymptotic expansion. Moreover, we presented a non-asymptotic calculation formula of the NML code length for exponential family models. This formula can be applied if we know the partition function. In addition, if we know the closed form of a distribution, the characteristic function of which is given by a function that can be expressed using the partition function, the calculation is reduced to one integral calculation, which is often analytically obtained.

\bibliography{ref}
\bibliographystyle{abbrv}

\end{document}